\documentstyle[12pt]{article}

\catcode`\@=11
\begin{document}
\centerline{\large\bf Two-parameter Quantum Affine Superalgebra $U_{p,q}\widehat{(gl(1|1))}$}
\centerline{\large\bf  and Its Drinfel$^{\prime}$d's Realization}
\vspace{0.8cm}
\centerline{\sf $^a$Jin-fang Cai\footnote{\sl  Email address: caijf@itp.ac.cn}
,  $^{bc}$Shi-kun Wang, $^a$Ke Wu }
\baselineskip=13pt
\vspace{0.5cm}
\centerline{$^a$ Institute of Theoretical Physics, Academia Sinica, }
\baselineskip=12pt
\centerline{ Beijing, 100080, P. R. China }
\vspace{0.3cm}
\baselineskip=13pt
\centerline{ $^b$ CCAST (World Laboratory), P.O. Box 3730, Beijing, 100080, 
  P. R. China  }
\baselineskip=12pt
\vspace{0.3cm}
\centerline{ $^c$ Institute of Applied Mathematics, Academia Sinica,}
\baselineskip=12pt
\centerline{Beijing, 100080, P. R. China }
\vspace{0.9cm}
\begin{abstract} {Using super RS construction method and Gauss decomposition,
we obtain Drinfel$^{\prime}$d's currents realization of two-parameter
quantum affine superalgebra  $U_{p,q}\widehat{(gl(1|1))}$ and get  co-product structure for these currents}
\end{abstract}
\baselineskip=14pt
\vspace{1cm}

For quantum affine algebras, there are three realizations or constructions:
Drinfel$^{\prime}$d-Jimbo realization \cite{DRI1,J2},  Drinfel$^{\prime}$d
\cite{DRI2} and RS \cite{RS} (or affine version of FRT \cite{FRT})
 realization. The explicit isomorphism between
Drinfel$^{\prime}$d and RS realization was established by Ding and Frenkel
\cite{DF} using Gauss decomposition. Recently, some attensions have been 
paid to the construction or definition of quantum affine superalgebras 
\cite{YA, C1, YZZ}. In this paper, we will use the super version of
 RS construction
method and Ding-Frenkel isomorphic map  to define a two-parameter quantum 
affine superalgebra  $U_{p,q}\widehat{(gl(1|1))}$ and get its
 Drinfel$^{\prime}$d's currents realization. We also calculate co-product 
structure for these currents.

 The graded Yang-Baxter equation(YBE) takes this form \cite{YH}: 
\begin{equation} 
\eta _{12}R_{12}(z/w) \eta _{13}R_{13}(z) \eta _{23}R_{23}(w)
= \eta _{23}R_{23}(w) \eta _{13}R_{13}(z) \eta _{12}R_{12}(z/w),
\label{ybe}\end{equation}
where $R$-matrix acts in tensor product of two 2-D graded linear space $V$:
$R(z) \in End(V \otimes V) $, and $\eta _{ik,jl}=(-1)^{(i-1)(k-1)} \delta _{ij} \delta _{lk}$ . It can also be writen in components as:
\begin{eqnarray}
&&R^{ij}_{ab}(z/w)R^{ak}_{pc}(z)R^{bc}_{qr}(w)
(-1)^{(P(a)-P(p))P(b)}   \nonumber    \\
&&{\hspace{2cm}}=(-1)^{P(e)(P(f)-P(r))}
R^{jk}_{ef}(w)R^{if}_{dr}(z)R^{de}_{pq}(z/w)
\end{eqnarray}

It's very obvious that $\eta R(z)$ satisfies the ordinary YBE when $R(z)$
is a solution of graded YBE. It can be verified that following $R$-matrix  
\cite{YH} is a solution of graded YBE (\ref{ybe}) :
\begin{equation} 
R_{12}(z)=\left(\begin{array}{lccr}
1&0 & 0& 0\\ 
0& \frac{(z-1)qp^{-1}}{zq-p^{-1}} & \frac{z(q-p^{-1})}{zq-p^{-1}} & 0\\
0& \frac{(q-p^{-1})}{zq-p^{-1}} & \frac{z-1}{zq-p^{-1}} & 0 \\
0&0&0&-\frac{q-zp^{-1}}{zq-p^{-1}} 
\end{array} \right)
\end{equation}
This solution satisfy the unitary condition: $R_{12}(z)R_{21}(z^{-1})=
{\bf 1}$, and it have two deformation parameters: $p$ and $q$. 
From the above solution  of graded YBE, we can define two-parameter quantum 
affine superalgebra  $U_{p,q}\widehat{(gl(1|1))}$ with a central extension 
employing the super version RS construction method \cite{RS}. 
 $U_{p,q}\widehat{(gl(1|1))}$ is an associative algebra with generators
\{$l_{ij}^k \vert 1\leq i,j \leq 2, k\in {\bf Z }$\} which subject to the 
following multiplication relations:
\begin{eqnarray}
&R_{12}(\frac{z}{w})L_1^{\pm}(z)  \eta L_2^{\pm}(w) \eta =
 \eta L_2^{\pm}(w) \eta L_1^{\pm}(z) R_{12}(\frac{z}{w})  \label{rll1}
\\
&R_{12}(\frac{z_-}{w_+})L_1^{+}(z)  \eta L_2^{-}(w) \eta =
 \eta L_2^{-}(w) \eta L_1^{+}(z) R_{12}(\frac{z_+}{w_-}) \label{rll2}
\end{eqnarray}
here $z_{\pm}=zq^{\pm \frac{c}{2}}$. We have used standard notation:
$L_1^{\pm}(z)=L^{\pm}(z)\otimes {\bf 1}, L_2^{\pm}(z)={\bf 1} \otimes
L^{\pm}(z) $ and $L^{\pm}(z)=\left( l_{ij}^{\pm}(z) \right)^2_{i,j=1}$,
 $l_{ij}^{\pm}(z)$ are generating functions (or currents) of $l_{ij}^k$:
$l_{ij}^{\pm}(z)=\sum_{k=0}^{\infty}l_{ij}^{\pm k} z^{\pm k}$.

This algebra admits coalgebra and antipole structure compatible with the 
associative \\
 multiplication defined by eqs.(4)(5): 
\begin{eqnarray}
&&\triangle \left(l_{ij}^{\pm}(z)\right) =\sum_{k=1}^2 l_{kj}^{\pm}
(zq^{\pm \frac{c_2}{2}})\otimes l_{ik}^{\pm}(zq^{\mp \frac{c_1}{2}})
(-1)^{(k+i)(k+j)},\label{cop}\\
&& \epsilon\left(l_{ij}^{\pm}(z)\right)=\delta _{ij} ,{\hspace{1.cm}}
S\left( ^{st} L^{\pm}(z) \right) =\left[ ^{st} L^{\pm}(z) \right]^{-1}.
\label{es}\end{eqnarray}
where  $\left[^{st} L^{\pm}(z) \right]_{ij} =(-1)^{i+j}l_{ji}^{\pm}(z) $.

  $L^{\pm}(z)$ have the following unique decompositions \cite{DF} : 
\begin{eqnarray}
L^{\pm}(z)& = & \left( \begin{array}{lr} 1& 0\\f^{\pm}(z) & 1 \end{array}
\right) \left( \begin{array}{lr}k_1^{\pm}(z) & 0\\0&k_2^{\pm}(z) \end{array}
\right) \left( \begin{array}{lr} 1 & e^{\pm}(z) \\ 0&1 \end{array} 
\right) 
\end{eqnarray}
where $ e^{\pm}(z), f^{\pm}(z) $ and $k_i^{\pm}(z)$ ($i$=1,2) is generating 
functions of  $U_{p,q}\widehat{(gl(1|1))}$ and $k_i^{\pm}(z)$ ($i$=1,2) are 
invertible. Let
$X^+(z)=e^+(z_-)-e^-(z_+) $ and $ X^-(z)=f^+(z_+)-f^-(z_-) $.
By the similar calculation made in \cite{DF, C1} for quantum affine (super)
algebras, we obtain the following 
(anti-)commutation relations among $ X^{\pm}(z)$ and $k_i^{\pm}(z)$ ($i$=1,2):
\begin{eqnarray}
&&[k_1^{\pm}(z) ~,~k_1^{\pm}(w)]=[k_1^+(z) ~,~k_1^-(w)]=0  \\
&&[k_1^{\pm}(z) ~,~k_2^{\pm}(w)]=[k_2^{\pm}(z) ~,~k_2^{\pm}(w)]=0  \\
&&\frac{z_{\pm}-w_{\mp}}{z_{\pm}q-w_{\mp}p^{-1}}k_1^{\pm}(z)
k_2^{\mp}(w)^{-1}=\frac{z_{\mp}-w_{\pm}}{z_{\mp}q-w_{\pm}p^{-1}}
k_2^{\mp}(w)^{-1}k_1^{\pm}(z)   \\
&&\frac{w_+q-z_-p^{-1}}{z_-q-w_+p^{-1}}k_2^+(z)^{-1}k_2^-(w)^{-1}
=\frac{w_-q-z_+p^{-1}}{z_+q-w_-p^{-1}}k_2^-(w)^{-1}k_2^+(z)^{-1} \\
&&k_i^{\pm}(z)^{-1}X^+(w)k_i^{\pm}(z)=\frac{z_{\pm}p-wq^{-1}}{
z_{\pm}-w}X^+(w) \hspace{1cm} (i=1,2)  \\
&&k_i^{\pm}(z)X^-(w)k_i^{\pm}(z)^{-1}=\frac{z_{\mp}p-wq^{-1}}{
z_{\mp}-w}X^-(w) \hspace{1cm} (i=1,2)  \\
&&\{X^+(z) ~,~X^+(w)\}=\{X^-(z) ~,~X^-(w)\}=0 \\
&&\{X^+(z) ~,~X^-(w)\}=(p-q^{-1})\left[\delta(\frac{w_-}{z_+})
k_1^-(z_+)^{-1}k_2^-(z_+)-\delta(\frac{z_-}{w_+})k_1^+(w_+)^{-1}
k_2^+(w_+)\right]
\end{eqnarray}
here $\delta (z)=\sum_{k\in {\bf Z}}z^k$.
The above relations are  Drinfel$^{\prime}$d's currents realization
 of two-parameter quantum affine superalgebra $U_{p,q}\widehat{(gl(1|1))}$.
 It's very clear that 
 $ X^{\pm}(z)$ (or $ e^{\pm}(z)$ and $f^{\pm}(z) $ ) are Fermionic type for
their anti-commutation relations and $k_i^{\pm}(z)$ ($i$=1,2) are Bosonic 
type elements in  $U_{p,q}\widehat{(gl(1|1))}$ as expected.

Introducing a transformation for currents:
\begin{eqnarray}
&&E(z)=X^+(zq) {\hspace{2.5cm}} F(z)=X^-(zq) \\
&&K^{\pm}(z)=k_1^{\pm}(zq)^{-1}k_2^{\pm}(zq) {\hspace{0.7cm}} 
H^{\pm}(z)=k_2^{\pm}(zq)k_1^{\pm}(zp^{-1})
\end{eqnarray}
then the (anti-)commutation relations for $E(z), F(z), K^{\pm}(z)$
and $H^{\pm}(z) $ are:
\begin{eqnarray}
&&[K^{\pm}(z) ~,~K^{\pm}(w)]=[H^{\pm}(z) ~,~ H^{\pm}(w)] =0 \\
&&[K^+(z) ~,~ K^-(w)] =[K^{\pm}(z) ~,~H^{\pm}(w)]=0 \\
&&\frac{w_{\pm}q-z_{\mp}p^{-1}}{w_{\pm}q^{-1}-z_{\mp}p}K^{\pm}(z)H^{\mp}(w)
=H^{\mp}(w)K^{\pm}(z)\frac{w_{\mp}q-z_{\pm}p^{-1}}{w_{\mp}q^{-1}-z_{\pm}p} \\
&&\left( \frac{z_+q-w_-p^{-1}}{z_+p^{-1}-w_-q}\right)^2 H^+(z)H^-(w)
=H^-(w)H^+(z)\left( \frac{z_-q-w_+p^{-1}}{z_-p^{-1}-w_+q}\right)^2  \\
&&[K^{\pm}(z) ~,~E(w)]=[K^{\pm}(z) ~,~F(w)]=0  \\
&&E(w)H^{\pm}(z)=\frac{z_{\pm}p-wq^{-1}}{z_{\pm}p^{-1}-wq}H^{\pm}(z)E(w) \\
&&H^{\pm}(z)F(w)=\frac{z_{\mp}p-wq^{-1}}{z_{\mp}p^{-1}-wq}F(w)H^{\pm}(z) \\
&&\{E(z) ~,~E(w)\}=\{F(z) ~,~F(w)\}=0 \\
&&\{E(z) ~,~F(w)\}=(p-q^{-1})\left[\delta(\frac{w_-}{z_+})
K^-(z_+)-\delta(\frac{z_-}{w_+})K^+(w_+)\right]
\end{eqnarray}
The co-product structure for currents  can be calculated directly 
from (\ref{cop})(\ref{es}):
\begin{eqnarray}
&&\triangle \left(K^{\pm}(z)\right)= K^{\pm}(zq^{\pm \frac{c_2}{2}}) 
\otimes  K^{\pm}(zq^{\mp \frac{c_1}{2}})  \\
&&\triangle \left(H^{\pm}(z)\right) = H^{\pm}(zq^{\pm \frac{c_2}{2} }) 
\otimes  H^{\pm}(zq^{\mp \frac{c_1}{2} }) -(p+q^{-1})F^{\pm}(zp^{-1}q^{-1\mp 
\frac{c_1}{2}\pm \frac{c_2}{2}}) H^{\pm}(zq^{\pm \frac{c_2}{2}})
   \otimes \nonumber \\
&&\hspace{3cm} H^{\pm}(zq^{\mp\frac{c_1}{2}})
E^{\pm}(zp^{-1}q^{-1\mp \frac{c_1}{2}\pm \frac{c_2}{2}})  \\
&&\triangle \left(E^{\pm}(z)\right)=E^{\pm}(z) \otimes {\bf 1}
+K^{\pm}(zq^{\mp\frac{c_1}{2}})\otimes E^{\pm}(zq^{\mp c_1})  \\
&&\triangle \left(F^{\pm}(z)\right)={\bf 1}\otimes F^{\pm}(z) 
+F^{\pm}(zq^{\pm c_2})\otimes K^{\pm}(zq^{\pm\frac{c_2}{2}}) \\
&&\epsilon \left(K^{\pm}(z)\right)=\epsilon \left(K^{\pm}(z)\right)=1 \\
&&\epsilon \left(E^{\pm}(z)\right)=\epsilon \left(F^{\pm}(z)\right)=0
\end{eqnarray}
where $E^{\pm}(z)=e^{\pm}(z_{\mp}q)$ and $F^{\pm}(z)= f^{\pm}(z_{\pm}q)$, then
$E(z)=E^+(z)-E^-(z)$ and $F(z)=F^+(z)-F^-(z)$.

\end{document}